\documentclass{amsart}
\usepackage{amsfonts, amsmath, latexsym, epsfig}
\usepackage[norelsize]{algorithm2e}
\usepackage{amssymb}
\usepackage{color}
\usepackage{epsf}
\usepackage{url}

\newcommand{\RR}{\ensuremath{\mathbb{R}}}

\newtheorem{proposition}{Proposition}
\newtheorem{theorem}{Theorem}

\newtheorem{conjecture}{Conjecture}

\newtheorem{remark}{Remark}
\newtheorem{definition}{Definition}
\def\QuotS#1#2{\leavevmode\kern-.0em\raise.2ex\hbox{$#1$}\kern-.1em/\kern-.1em\lower.25ex\hbox{$#2$}}

\DeclareMathOperator{\Aut}{Aut}

\DeclareMathOperator{\HYP}{HYP}
\DeclareMathOperator{\CUT}{CUT}

\DeclareMathOperator{\OMCUTP}{OMCUTP}
\DeclareMathOperator{\CUTP}{CUTP}

\DeclareMathOperator{\MET}{MET}
\DeclareMathOperator{\METP}{METP}
\DeclareMathOperator{\QMET}{QMET}
\DeclareMathOperator{\QMETP}{QMETP}
\DeclareMathOperator{\WQMET}{WQMET}
\DeclareMathOperator{\WQMETP}{WQMETP}
\DeclareMathOperator{\WOMCUTP}{WOMCUTP}
\DeclareMathOperator{\OMCUT}{OMCUT}
\DeclareMathOperator{\WOMCUT}{WOMCUT}
\DeclareMathOperator{\HMET}{HMET}

\begin{document}

\author{Michel Deza}
\address{Michel Deza, \'Ecole Normale Sup\'erieure, Paris, Deceased}

\author{Mathieu Dutour Sikiri\'c}
\address{Mathieu Dutour Sikiri\'c, Rudjer Boskovi\'c Institute, Bijenicka 54, 10000 Zagreb, Croatia, Fax: +385-1-468-0245}
\email{mathieu.dutour@gmail.com}

\title{Generalized cut and metric polytopes of graphs and simplicial complexes}
\date{}

\begin{abstract}
  Given a graph $G$ one can define the cut polytope $\CUTP(G)$ and the metric polytope $\METP(G)$ of this graph
  and those polytopes encode in a nice way the metric on the graph.
  According to Seymour's theorem, $\CUTP(G) = \METP(G)$ if and only if $K_5$ is not a minor of $G$.

  We consider possibly extensions of this framework:
  \begin{enumerate}
  \item We compute the $\CUTP(G)$ and $\METP(G)$ for many graphs.
  \item We define the oriented cut polytope $\WOMCUTP(G)$ and oriented multicut polytope $\OMCUTP(G)$ as well as their oriented metric version $\QMETP(G)$ and $\WQMETP(G)$.
  \item We define an $n$-dimensional generalization of metric on simplicial complexes.
  \end{enumerate}

\end{abstract}


\maketitle

\section{Introduction}\label{Sec_Introduction}
The cut polytope \cite{DL} is a natural polytope arising in the study of the maximum cut problem \cite{CompendiumOptiPb}.
The cut polytope on the complete graph $K_n$ has seen much study (see \cite{DL}) but the cut polytope on
a graph was much less studied \cite{CUTsmallGraphs,BarahonaCutContract,LeggettGargIneq}.
Moreover, generalizations of the cut polytope on graphs seems not to have been considered.

Given a graph $G=(V,E)$, for a vertex subset $S\subseteq V=\{1, \dots, n\}$, the {\em cut semimetric} $\delta_S(G)$
is a vector (actually, a symmetric $\{0,1\}$-matrix) defined as
\begin{equation}\label{DefCutMetric}
\delta_S(x,y)=\left\{\begin{array}{rl}
1 &\mbox{if~~~} (xy)\in E \mbox{~~~~and ~~~}  \vert S \cap \{x,y\}\vert = 1,\\
0 &\mbox{otherwise.}
\end{array}\right.
\end{equation}
A {\em cut polytope} $\CUTP(G)$, respectively {\em cut cone} $\CUT(G)$, are defined as the convex hull of all such semimetrics, respectively positive span of all non-zero ones among them.
The dimension of $\CUTP(G)$ and $\CUT(G)$ is equal to the number of edges of $G$.

The {\em metric cone} $\MET(K_n)$ is the set of all {\em semimetrics} on $n$ points, i.e.,
the functions $d : \{1,\dots , n\}^2 \rightarrow \mathbb{R}_{\ge 0}$
(actually, symmetric matrices over $\mathbb{R}_{\ge 0}$ having only zeroes on the diagonal), 
which satisfy  all $3 {n\choose 3}$
{\em triangle inequalities}  $d(i,j) + d(i,k) - d(j,k) \ge 0$.
 The bounding of $\MET(K_n)$ by $n\choose 3$ {\em perimeter inequalities}
$d(i,j) + d(i,k) + d(j,k) \le 2$ produces the {\em metric polytope} $\METP(K_n)$.

For a graph $G=(V,E)$ of the {\em order} $|V|=n$, let $\MET(G)$ and $\METP(G)$ denote the projections of 
$\MET(K_n)$ and $\METP(K_n)$, respectively, on the subspace $\mathbb{R}^{|E|}$ indexed by the edge set of $G$.
Clearly, $\CUT(G)$ and $\CUTP(G)$ are projections of, respectively, $\CUT(K_n)$ and $\CUTP(K_n)$
on $\mathbb{R}^E$.
It holds
\begin{equation*}
\CUT(G)\subseteq \MET(G)  \mbox{~and~} \CUTP(G)\subseteq \METP(G).
\end{equation*}

In Section \ref{Section_CutPolytope} we consider the structure of those polytopes and give
the description of the facets for many graphs (see Tables \ref{TableSomeCUTpolytopes1} and
\ref{TableSomeCUTpolytopes2}).
The data file of the groups and orbits of facets of considered polytopes is available
from \cite{WebPageCutPolytopes}.

The construction of cuts and metrics can be generalized to metrics which are not necessarily
symmetric are considered in Section \ref{Sec_Quasi} (see also \cite{QuasiMetric1,QuasiMetric2}).
The triangle inequality becomes $d(i,j) \leq d(i,k) + d(k,j)$ and the perimeter
inequality becomes $d(i,j) + d(j,k) + d(k,i) \leq 2$ for $1\leq i,j,k\leq n$. We also need
the inequalities $0\leq d(i,j) \leq 1$.
The quasi metric polytope $\QMETP(K_n)$ is defined by the above inequalities and the quasi metric
cone $\QMET(K_n)$ is defined by the inequalities passing by zero.
The quasi metric cone $\QMET(G)$ and polytope $\QMETP(G)$ are defined as projection of above two
cone and polytopes. In Theorem \ref{QMET_theorem} we give an inequality description of those
projections.

Given an ordered partition $(S_1, \dots, S_r)$ of $\{1, \dots, n\}$ we defined an oriented
multicut as:
\begin{equation*}
  \delta'(S_1, \dots, S_r)_{x,y} = \left\{\begin{array}{rl}
  1 & \mbox{~if~}  x\in S_i, y\in S_j \mbox{~and~} i<j,\\
  0 & \mbox{~otherwise}.
  \end{array}\right.
\end{equation*}
The convex cone of the oriented multicut is the oriented multicut cone $\OMCUT(K_n)$.
The convex polytope can also be defined but there are vertices besides the oriented multicuts.
A smaller dimensional cone $\WQMET(G)$ and polytope $\WQMETP(G)$ can be defined by adding the cycle equality
\begin{equation*}
d(i,j) + d(j,k) + d(k,i) = d(j,i) + d(k,j) + d(i,k)
\end{equation*}
to the cone $\QMET(G)$ and polytope $\QMETP(G)$.
A multicut satisfies the cycle equality if and only if $r=2$. We note the corresponding cone $\WOMCUT(G)$
and $\WOMCUTP(G)$. In Section \ref{Sec_Quasi} we consider those cones and polytopes and their facet
description.

The notion of metrics can be generalized to more than $2$ points and we obtain the {\em hemimetrics}.
Those were considered in \cite{DezaRosenberg1,DezaRosenberg2,DezaDutourHemimetric1,DD_mining}. Only
the notion of cones makes sense in that context. The definition of the above papers extends the triangle
inequality in a direct way: It becomes a simplex inequality with the area of one side being bounded by
the sums of area of the other sides.
In \cite{BookGenFiniteMetrCuts} we argued that this definition was actually inadequate since it prevented
right definition of hemimetric for simplicial complex. In Section \ref{Sec_Simplicial} we give full
details on what we argue is the right definition of hemimetric cone.

There is much more to be done in the fields of metric cones on graphs and simplicial complexes.
Besides further studies of the existing cones and the ones defined in this paper, two other cases
could be interesting.
One is to extend the notion of hypermetrics cone $\HYP(K_n)$ to graphs; several approaches were considered
in \cite{GeneralizationHypermetric}, for example projecting only on the relevant coordinates, 
but no general results were proved.

Another generalization that could be considered is the {\em diversities} considered
in \cite{BT_hyperconvexity,BT_hyperconvexityOTHER}. {\em Diversity cone $DIV_n$} is the
set of all {\em diversities on $n$ points}, i.e., the functions
$f: \{A: A \subseteq  \{1,\dots,n\}\} \rightarrow \RR_{\ge 0}$ satisfying $f(A)= 0$ if $|A|\le 1$ and
\begin{equation*}
f(A\cup B)+f(B\cup C)\ge f(A\cup C) \,\mbox{~if~} \,B\neq \emptyset.
\end{equation*} 
The {\em  induced diversity metric} $d(i,j)$ is $f(\{i,j\})$.

{\em Cut diversity cone $CDIV_n$} is the positive span of all {\em cut diversities} $\delta (A)$,
where $A \subseteq  \{1,\dots,n\}$, which are defined, for any $S \subseteq  \{1,\dots,n\}$,  by  
\begin{equation*}
\delta_{S}(A) = \left\{\begin{array}{rl}
1 & \mbox{if~} A\cap S \neq \emptyset \mbox{~and~} A\setminus S \neq \emptyset,\\
0 & \mbox{otherwise.}
\end{array}\right.
\end{equation*} 

$CDIV_n$ is the set of all diversities from  $DIV_n$, which are isometrically embeddable
into an {\em $l_1$-diversity}, i.e., one, defined on $\mathbb{R}^m$
with $m \le {n\choose \left\lfloor \frac{n}{2} \right\rfloor}$ by
\begin{equation*}
f_{m1}(A)= \sum_{i=1}^m \max_{a,b\in A}\{|a_i-b_i|\}.
\end{equation*} 

These two cones are extensions of the $\MET(K_n)$ and $\CUT(K_n)$ on a complete hypergraphs
and it would be nice to have a nice definition on any hypergraph.

\section{Structure of cut polytopes of graphs}\label{Section_CutPolytope}

The cut metric $\delta_S$ defined at Equation \eqref{DefCutMetric} satisfies the relation
$\delta_{\{1, \dots, n\} - S} = \delta_S$. The cut polytope $\CUTP(K_n)$ is defined
as the convex hull of the metrics $\delta_S$ and thus has $2^{n-1}$ vertices.

For a given subset $S$ of $\{1, \dots, n\}$ we can define the {\em switching} operation $F_S$
by
\begin{equation*}
  F_S(d)(i,j) = \left\{\begin{array}{cl}
  1 - d(i,j) & \mbox{if~} \left\vert S\cap \{i,j\}\right\vert = 1,\\
  d(i,j)     & \mbox{otherwise}.
  \end{array}\right.
\end{equation*}
The operation on cuts is $F_S(\delta_T) = \delta_{S\Delta T}$ with $\Delta$ denoting the symmetric difference
(see \cite{DL} for more details).
For a graph $G$ we define $\CUTP(G)$ to be the projection of $\CUTP(K_n)$ on the
coordinates corresponding to the edges of the graph $G$. If $G$ is connected then $\CUTP(G)$ has
exactly $2^{n-1}$ vertices.
Then $\delta_S$ can be seen also as the adjacency matrix of a {\em cut}  (into $S$ and $\overline{S}$) {\em subgraph} of $G$.
The cut cone $\CUT(G)$ is defined by taking the convex cone generated by the metrics $\delta_S$ but it is generally
not used in that section.

In fact, $\CUT(K_n)$ is the set of all $n$-vertex semimetrics, which embed isometrically into some metric space $l_1$, and rational-valued elements of $\CUT(K_n)$ correspond exactly to the $n$-vertex semimetrics, which embed isometrically, {\em up to a scale} $\lambda \in \mathbb{N}$, into 
the path metric of some $m$-cube $K_2^m$. It shows importance of this cone in Analysis and Combinatorics. 
The enumeration of orbits of facets of $\CUT(K_n)$ and $\CUTP(K_n)$ for $n\le 7$ was done in \cite{OrbitFacetCutPolytope5,OrbitFacetCutPolytope6,OrbitFacetCutPolytope7} for $n=5$, $6$, $7$ respectively, and in \cite{CR}, completed by \cite{CUTsmallGraphs}, for $n=8$.

\subsection{Automorphism group of cut polytopes}

The symmetry group $\Aut(G)$ of a graph $G=(V,E)$ induces symmetry of $\CUTP(G)$. For any $U\subset \{1,\dots, n\}$, the map $\delta_S\mapsto \delta_{U\Delta S}$ also defines a symmetry of $\CUTP(G)$. Together, those form the {\em restricted symmetry group} $ARes(\CUTP(G))$ of order $2^{|V|-1}\vert \Aut(G)\vert$. 
The full symmetry group $\Aut(\CUTP(G))$ may be larger.
In Tables \ref{TableSomeCUTpolytopes1},  \ref{TableSomeCUTpolytopes2}, such cases are marked by ${}^*$.
 Denote $2^{1-|V|}|\Aut(\CUTP(G))|$ by $A(G)$.
 
For example, $|\Aut(\CUTP(K_n))|$ is $2^{n-1}n!$ if $n\neq 4$ and $6\times 2^34!$ if $n=4$ (\cite{Relatives}).

\begin{remark}

(i) If $G=(V,E)$ is $Prism_m$ ($m\neq 4$), $APrism_m$ ($m> 3$),  M\"{o}bius ladder $M_{2m}$ and $Pyr^2(C_m)$ ($m> 3$), then $Aut(G)=4m$.

(ii) If $G$ is a complete multipartite graph with $t_1$ parts of size $a_1$, $\dots ,$ $t_r$ parts of size $a_r$,
with $a_1< a_2< \dots <a_r$ and all $t_i\ge 1$, then $|Aut(G)|=\prod_{i=1}^rt_i!(a_i!)^{t_i}.$

(iii) Among the  cases considered here,
 all occurrences of $A(G)> |Aut(G)|$ are:
 $A(G)=m!2^{m-1}|Aut(G)|$ for $G=K_{2,m>2},K_{1,1,m>1}$ 
and $A(G)=6|Aut(G)|$, i.e., $2m!=48, 6m!$ for $G=K_{2,2}$ and  $K_{1,1,1,1}$, respectively.

(iv) 
If $G=P_m$ ($m\ge 3$ edges), then $|Aut(G)=2$, while $A(G)=m!=(|V|-1)!$. 

If $G=C_m$ ($m>3$), then $|Aut(G)|=2m$, while $A(G)=2m!$ for $m=4$ and $A(G)=m!=|V|!$ for $m\ge 5$.


\end{remark}

\subsection{Edge faces, $s$-cycle faces and metric polytope}

\begin{definition}
Let $G=(V,E)$ be a graph.

(i) Given an edge $e\in E$, the {\em edge inequality} (or
 {\em $2$-cycle inequality}) is
\begin{equation*}
x(e) \geq 0.
\end{equation*}

(ii) Given a $s$-cycle $c=(v_1, \dots, v_s), s\ge 3,$ of $G$, the {\em $s$-cycle inequality} is:
\begin{equation*}
x(c, (v_1, v_s)) = \sum_{i=1}^{s-1} x(v_i,v_{i+1})  -x(v_1,v_s)\ge 0.
\end{equation*}

\end{definition}
The edge inequalities and $s$-cycle inequalities are valid  on $\CUTP(G)$, since they are, clearly, valid on each cut:
a cut  intersects a cycle in the set of even cardinality.
So, they define faces, but not necessarily facets. In fact, it holds

\begin{theorem}

(i) The inequality $x(e)$ is facet defining in $\CUTP(G)$ (also, in $\CUT(G)$)
if and only if $e$ is not contained into a $3$-cycle of $G$.
 
(ii) An $s$-cycle inequality is facet defining in $\CUTP(G)$ (also, in $\CUT(G)$)
if and only corresponding $s$-cycle is chordless.

(iii) $\METP(G)$ is defined by all edge and $s$-cycle inequalities, while
$\MET(G)$ is defined by all  $s$-cycle inequalities.

\end{theorem}

In fact, (i) and (ii) above were proved in \cite{BarahonaMahjoubOnCutPolytope},
(iii) was proved in \cite{BarahonaCutsMatchingsPlanar}; see also Section 27.3 in \cite{DL}.

The following Theorem, proved in \cite{SeymourMatroidMulticommodity} for cones
and in \cite{BarahonaCutContract} for polytopes, clarifies when the metric and
cut polytope coincides:

\begin{theorem}\label{K5wonder}
$\CUT(G) = \MET(G)$ or, equivalently,  $\CUTP(G)=\METP(G)$ if and only if $G$ does not have any $K_5$-minor.
\end{theorem}

As a corollary of Theorem \ref{K5wonder}, we have that the facets of $\CUTP(G)$ (also, in $\CUT(G)$)
are determined by edge inequalities and $s$-cycle inequalities if and only if $G$  does not have any $K_5$-minor.

$3$-cycle inequality is usual triangle inequality; in fact, it is unique, among edge and all $s$-cycle inequalities
to define a facet in a $\CUTP(K_n)$.
 
The {\em girth} and {\em circumference} of a graph, having cycles, are the length of the shortest and longest cycle, respectively.
In a graph  $G$, a {\em chordless cycle} is any cycle, which is induced subgraph; so, any triangle, any shortest cycle and any cycle, bounding a face in some embedding of $G$, are chordless.
 Let $c_s'$ and $c_s$ denote the number of all and of all chordless  $s$-cycles in $G$, respectively.

There are $2|E|$ edge faces, which decompose into orbits, one for each orbit of edges  of $G$ under $\Aut(G)$.
 There are $2^{s-1}c'_s$ $s$-cycle faces, which decompose into orbits, one for each orbit of  $s$-cycles  of $G$ under $\Aut(G)$.

The incidence of  edge faces is $2^{|V|-2}$ and the size of each orbit is twice the size of corresponding orbit  of edges.
The incidence of  $s$-cycle faces is $2^{|V|-s}s$ and the size of each orbit is $2^{s-1}$ times
 the size of corresponding orbit  of $s$-cycles in $G$.

By {\em Wagner's theorem} \cite{Wagner1937}, a finite graph is planar if and only if it has no minors
$K_5$ and $K_{3,3}$. For embeddability on the projective plane $\mathbb{P}^2$, there are
exactly $103$ forbidden topological minors and exactly $35$ forbidden minors
(see \cite{ArchDeacon35,GloverHuneke103}).
For embeddability on the torus $\mathbb{T}^2$, $16629$ forbidden minors are known
(see \cite{GagarinMyrvoldChambers}) but the list is not necessarily complete.
Closely related {\em Kuratowski's theorem}
\cite{Kuratowski1930} states that a finite graph is planar if and only if it does not contain
a subgraph that is a subdivision of $K_5$ or of $K_{3,3}$.

\begin{table}
\begin{center}
\caption{The number of facets of  $\CUTP(G)$ of some $K_5$-minor-free graphs $G$; ${}^{*}$ shows $A(G)$=$2^{1-|V|} |\Aut(\CUTP(G))|>|Aut(G)|$}
\label{TableSomeCUTpolytopes1}
\begin{tabular}{||c|c||c|c|c||}
\hline\hline
 $G=(V,E)$       & $\vert V\vert,\vert E\vert$ & $A(G)$ & Number of facets&Orbit's $s$\\\hline\hline
M\"{o}bius ladder $M_{8}$& $8,2$              & $16$   & $184(4)$&$2,2,4,5$\\
$M_{6}=K_{3,3}$& $6,9$              & $2(3!)^2$  & $90(2)$&$2,4$\\\hline
$K_{1,1,1,m}, m> 1$     & $m$+$3,3m$+$3$ & $3!m!$ &$4+12m(2)$&$3,3$\\
$K_{1,2,m},m>1$   & $m$+$3,3m$+$2$          &  $|Aut(K_{1,2,m})|$  & $8m + 8{m\choose 2}(2)$&$3,4$\\
$K_{3,m}, m\ge 3$  & $m$+$3,3m$          & $|Aut(K_{3,m})|$      &  $6m + 24{m\choose 2}(2)$&$2,4$\\
$K_{2,m}, m> 2$ & $m$+ $2,2m$          & $2^{m-1}m!|Aut(K_{2,m})|$ $\,^{*}$  & $4m^2(1)$&$2$ with $4$\\
$K_{2,2}$    & $4,4$          & $6|\Aut(K_{2,2})|$  $\,^{*}$  & $16(1)$&$2$ with $4$\\ 
$K_{1,1,m}, m>1$     & $m$+$2,2m$+$1$ & $2^{m-1}m!|Aut(K_{1,1,m})|$ $\,^{*}$ &$4m(1)$ & $3$\\
$K_{m+1}$-$K_m$=$K_{1,m}, m$$>$$1$     & $m$+$1,m$ &  $m!$     &$2m(1)$ &$2$\\
\hline
$APrism_6$     & $12,24$          & $24$    & $2,032(5)$&$3,6,7,7,8$\\
$APrism_5$     & $10,20$          & $20$    & $552(4)$&$3,5,6,7$\\
$APrism_4$     & $8,16$          & $16$     & $176(3)$&$3,4,5$\\
$Prism_7$     & $14,21$          & $28$     & $7,394(6)$&$2,2,4,7,9,9$\\
$Prism_6$     & $12,18$          & $24$     & $2,452(6)$&$2,2,4,6,8,8$\\
$Prism_5$     & $10,15$          & $20$  & $742(5)$&$2,2,4,5,7$\\
$Prism_3$     & $6,9$          & $12$   & $38(3)$&$2,3,4$\\\hline
Tr. Tetrahedron  & $12,18$              & $24$    & $540(4)$&$2,3,6,8$\\
Cuboctahedron  & $12,24$              & $48$    & $1,360(5)$&$3,4,6,6,8$\\
\hline
Dodecahedron  & $20,30$            & $120$    & $23,804(5) $&$2,5,9,10,10$\\
Icosahedron  & $12,30$            & $120$    & $1,552(4)$&$3,5,6,6$\\
Cube   $K_2^2$     & $8,12$              & $48$    & $200(3)$&$2,4,6$\\
Octahedron $K_{2,2,2}$ & $6,12$              & $48$    & $56(2)$&$3,4$\\
Tetrahedron $K_4$ & $4,6$              & $6|Aut(K_4)|$\,\,\,${}^{*}$    & $12(1)$&$3$\\
\hline\hline
\end{tabular}
\end{center}
\end{table}

\begin{table}
\begin{center}
\caption{The number of facets of  $\CUTP(G)$ for some graphs $G$ with $K_5$-minor}
\label{TableSomeCUTpolytopes2}
\begin{tabular}{||c|c||c|c|c||}
\hline\hline
 $G=(V,E)$       & $\vert V\vert,\vert E\vert$ & $A(G)$ & Number of facets (orbits)&Orbit's $s$\\\hline\hline
Heawood graph &$14,21$         & $336$ & $5,361,194(9)$&$2,6,8$\\
Petersen graph & $10,15$        & $120$  & $3,614(4)$&$2,5,6$\\\hline
M\"{o}bius ladder $M_{10}$& $10,15$              & $20$   & $1,414(5)$&$2,2,4,6$\\\hline
M\"{o}bius ladder $M_{12}$ & $12,18$              & $24$   & $26,452(6)$&$2,2,4,7,9$\\
M\"{o}bius ladder $M_{14}$ & $14,21$              &  $28$  & $369,506(9)$&$2,2,4,8,10$\\
$K_{5,5}$& $10,25$              & $2(5!)^2$  & $16,482,678,610(1,282)$&$2,4$\\
$K_{4,7}$& $11,28$              &  $4!7!$    & $271,596,584(15)$&$2,4$\\
$K_{4,6}$& $10,24$              &$4!6!$     & $23,179,008(12)$&$2,4$\\
$K_{4,5}$& $9,20$              & $4! 5!$    & $983,560(8)$&$2,4$\\
$K_{4,4}$& $8,16$              & $2(4!)^2$  & $27,968(4)$&$2,4$\\\hline
$K_{3,3,3}$&$9,27$&$(3!)^4$    &$624, 406, 788(2,015)$&$3,4$\\
$K_{1,4,4}$&$9,24$&$2(4!)^2$   &$36,391,264(175)$&$3,4$\\
$K_{1,3,5}$&$9,23$&$3!5!$      &$71,340(7)$&$3,4$\\
$K_{1,3,4}$&$8,19$&$3!4!$      &$12,480(6)$&$3,4$\\
$K_{1,3,3}$&$7,15$&$2(3!)^2$    &$684(3)$&$3,4$\\\hline
$K_{1,1,3,3}$&$8,21$&$4(3!)^2$  &$432,552(50)$&$3,3,4$\\
$K_{1,2,2,2}$&$7,14$&$3!(2!)^3$ &$5,864(9)$&$3,3,4$\\
$K_{1,1,2,2}$ &$6,13$&$4(2!)^2$ &$184(4)$&$3,3,4$\\
$K_{1,1,2,m},m> 2$&$m$+$4,4m$+$5$ & $4m!$   &$8$+$20m$+$8{m\choose 2}(16m-15)(7)$&$3,3,3,4$\\\hline
$K_{1,1,1,1,m}, m>1$ & $m$+$4,4m$+$6$ & $4!m!$   &$8(8m^2-3m+2)$(4)&$3,3$\\\hline
$K_{1,1,1,1,1,3}$=$K_8-K_3$& $8.25$              & $360$   & $2,685,152(82)$&$3,3$\\
$K_{1,1,1,1,1,2}$=$K_7-K_2$& $7,20$              & $240$  & $31,400(17)$&$3,3$\\
$K_7-C_3$& $7,18$              & $144$  & $520(4)$&$3,3$\\
$K_7-C_4$&$7,17$&$48$&$108$(4)&$3,3,3$\\ 
$K_7-C_5$=$Pyr^2(C_5)$&$7,16$&$20$&$780(6)$&$3,3,5$\\ 
$K_7-C_6$=$Pyr(Prism_3)$&$7,15$&$12$&$452(5)$&$3,3,3,4$\\ 
$K_7-C_7$&$7,14$&$14$&$148(3)$&$3,4$\\ \hline
$Pyr(Prism_4)$&$9,20$&$48$&$10,464(6)$&$3,4,6$\\
$Pyr(Prism_5)$&$11,25$&$20$&$208,133(22)$&$3,3,4,5,7$\\
$Pyr(APrism_4)$&$9,24$&$16$&$389,104(17)$&$3,3,3,4,5$\\\hline
$Pyr^2(C_6)$&$8,19$&$24$&$3,432(7)$&$3,3,6$\\
$Pyr^2(C_7)$&$9,22$&$28$&$14,740(11)$&$3,3,7$\\\hline
Tr.Octahedron on $\mathbb{P}^2$&$12,18$&$48$&$62,140(7)$&$2,2,4,6,6$\\
\hline\hline
\end{tabular}
\end{center}
\end{table}



\subsection{Skeletons of Platonic and semiregular polyhedra}

Let $G$ be embedded in some oriented surface; so, it is a map $(V,E,F)$, where $F$ is the  set of faces of $G$. Let $\vec{p}=(\dots, p_i, \dots)$ denote the {\em $p$-vector} of the map, enumerating the number $p_i>0$ of faces of all sizes $i$, existing in $G$.  

 Call {\em face-bounding} any $s$-cycle of $G$, bounding a face in map $G$.
Call an $s$-cycle of $G$ {\em $i$-face-containing},  {\em edge-containing} or
 {\em point-containing}, 
 if all its interior points form just $i$-gonal face, edge or   point,
 respectively. 
   Call {\em equator} any   cycle $C$, the interior of which (plus $C$) is isomorphic to the 
exterior (plus $C$).

The chordless $4,6,5,9$-cycles of Octahedron, Cube, Icosahedron and Dodecahedron, respectively,
are exactly their vertex-containing $4,6,5,9$-cycles.
 
For Octahedron and Cube, they are exactly all $3$ and $4$ equators, respectively, which are, apropos,  the {\em central circuits} and {\em zigzags}  (see \cite{newBook}), respectively.
 
All $c_6$ chordless  $6$-cycles of  Icosahedron  are exactly
their $30$ edge-containing ones and $10$ face-containing ones, which
are exactly the $10$ equators and the {\em weak zigzags} (\cite{newBook}).
All $c_{10}$ chordless  $10$-cycles of  Dodecahedron are $30$ edge-containing
ones and $6$ face-containing ones, which are exactly all $6$ equators and the
zigzags.

\begin{proposition}
If $G$ is the skeleton of a Platonic solid, then all possible facets of $\CUTP(G)$ are:
edge facets and  $s$-cycle facets, coming from  all  face-bounding cycles and from all
(if they exist and not listed before) vertex-, edge-, face-containing cycles.

For instance:

(i) If  $G=K_4$ (Tetrahedron), then $\CUTP(G)$ has 
unique orbit of $2^2p_3=16$ (simplicial) $3$-cycle facets (from all $|F|=p_3=4$ face-bounding cycles of $G$).

(ii) If $G=K_{2,2,2}$ (Octahedron), then $\CUTP(G)$ has $56$ facets in $2$ orbits, namely:   

orbit of  $2^2p_3$  $3$-cycle facets (from  all $|F|=p_3=8$ face-bounding cycles, 
orbit of $2^3c_4$  $4$-cycle facets (from  all $c_4 =\frac{|V|}{2}=3$ vertex-containing $4$-cycles).
 
(iii) If $G=K_{2}^3$ (Cube), then $\CUTP(G)$ has $200$ facets in $3$ orbits, namely:
orbit of  $2|E|=24$ edge facets, 
orbit of  $2^3p_4$  $4$-cycle facets (from  all $|F|=p_4=6$ face-bounding cycles), 
   
orbit of $2^5c_6=128$  $6$-cycle facets (from  all $c_6=4$ vertex-containing $6$-cycles).

(iv) If $G$ is Icosahedron, then $\CUTP(G)$ has $1,552$ facets in $4$ orbits, namely:

orbit of  $2^2p_3=80$  $3$-cycle facets (from  all $|F|=p_3=20$ face-bounding cycles), 
   
orbit of $2^4c_5=192$  $5$-cycle facets (from  all $c_5=12$ vertex-containing $5$-cycles), 
  
orbit of $2^5|E|=960$ $6$-cycle facets (from $|E|=30$  edge-containing $6$-cycles), 
  
orbit of $320$ $6$-cycle facets (from $\frac{|F|}{2}=10$ face-containing $6$-cycles).
   
(v) If $G$ is Dodecahedron, then $\CUTP(G)$ has $23,804$ facets in $5$ orbits, namely:

orbit of $2|E|=60$ edge facets, 
orbit of $2^4p_5=192$  $5$-cycle facets (from  all $|F|=p_5=12$ face-bounding cycles), 
orbit of $2^8c_9=5,120$ $9$-cycle facets (from  all $c_9=20$ vertex-containing $9$-cycles),  
orbit of $2^9|E|=15,360$ $10$-cycle facets (from $30$  edge-containing $10$-cycles), 
orbit of $2^9\times 6=3,072$ $10$-cycle facets (from $\frac{|F|}{2}=6$ face-containing $10$-cycles).
\end{proposition}

\bigskip

In a Truncated Tetrahedron, call {\em ring-edges} those bounding a triangle, and {\em rung-edges}
all $6$ other ones.

\begin{proposition}
(i) If $G$ is Truncated Tetrahedron, then $\CUTP(G)$  has $540$ facets:
\begin{enumerate}
\item orbit of  $2\times 6$ edge facets (from all $6$ rung-edges),
 
\item orbit of $2^3p_3$ $3$-cycle 
 facets (from  all $p_3=4$ $3$-face-bounding cycles), 
 
\item orbit of $2^5p_6$ $6$-cycle 
 facets (from  all $p_6=4$ $6$-face-bounding cycles), 
 
 \item  orbit of $2^7\times 3$
 $8$-cycle facets (from $\frac{1}{2}{4\choose 2}$ rung-edge-containing $8$-cycles,  which are also the equators).
 \end{enumerate} 
 (ii) If $G$ is Cuboctahedron, then $\CUTP(G)$ has $1,360$ facets, namely:
\begin{enumerate} 
\item orbit of $2^2p_3$ $3$-cycle 
 facets (from  all $p_3=8$  $3$-face-bounding cycles), 
 
\item orbit of $2^3p_4$ $4$-cycle facets (from all $p_4=6$ $4$-face-bounding cycles), 
 
\item orbit of $2^5|V|$ $6$-cycle facets (from all $12$ vertex-containing $6$-cycles), 
  
\item orbit of $2^5\times 4=128$  $6$-cycle facets (from all $\frac{p_3}{2}$ $3$-face-containing $6$-cycles, which are also equators and the central circuits),

\item orbit of $2^7p_4=768$  $8$-cycle facets (from all $6$ $4$-face-containing $8$-cycles, which are also zigzags).
\end{enumerate}

\end{proposition}

\bigskip

Given a $Prism_m$ ($m\neq 4$)  or an $APrism_m$ ($m\neq 3$), we  
 call {\em rung-edges} the edges  connecting two $m$-gons, and {\em ring-edges} other $2m$ edges. 
 
 Let $P$ be an ordered partition $X_1\cup \dots \cup X_{2t}=\{1, \dots , m\}$
 into ordered sets $X_i$ of $|X_i|\ge 3$ consecutive integers. Call
{\em $P$-cycle of $Prism_m$} the chordless $(m+2t)$-cycle obtained by taking
 the path $X_1$ on the, say, $1$-st $m$-gon, then rung edge (in the same
direction, then path $X_2$ on the $2$-nd $m$-gon, etc. till returning to
the path $X_1$. Any vertex of $Prism_m$ can be taken as the $1$-st element
of $X_1$, in order to fix a $P$-cycle. So, a $P$-cycle defines an orbit of
$2^{m+2t-1}2m$ $(m+2t)$-cycle facets of $\CUTP(Prism_m)$, except the case
$(|X_1|,\dots, |X_m|)=(|X_2|,\dots , |X_{2t}|,|X_1|)$ when the orbit is twice smaller.

A {\em $P$-cycle of $APrism_m$} is defined similarly, but we ask only
 $|X_i|\ge 2$ and rung edges, needed to change $m$-gon, should be
selected, in the cases $|X_i|= 2,3$ so that they not lead to a ring
edge,i.e., a chord
on $P$. Clearly, $P$-cycles are are all possible chordless $t$-cycles 
with $t\neq 4,m$ for   $Prism_m$ and with $t\neq 2,m$ for   $APrism_m$.

\begin{proposition} 
  
(i)   If $G$ is $Prism_m$ ($m\ge 5$),  then    all facets of $\CUTP(G))$ are:
\begin{enumerate}
\item orbit of $2m$ edge facets (from all $m$ rung-edges)
 \item orbit of $4m$ edge facets (from all $2m$ ring-edges);
 \item orbit of $2^3p_4=8m$  $4$-cycle facets (from all $m$  $4$-face-bounding $4$-cycles);
\item orbit of $ 2^{m-1}p_m$ of $m$-cycle facets (from both $m$-face-bounding $m$-cycles);
\item orbits of cycle facets for all possible $P$-cycles.

\end{enumerate}


(ii) If $G$ is $APrism_m (m\ge 4$),  then all facets of   $\CUTP(G))$ are:
 
\begin{enumerate}

 \item orbit of $2^2p_3=8m$  $3$-cycle facets (from all $2m$  $3$-face-bounding  $3$-cycles);
\item orbit of $2^{m-1}p_m$ of $m$-cycle facets (from both $m$-face-bounding $m$-cycles);

\item orbits of cycle facets for all possible $P$-cycles.

\end{enumerate}

\end{proposition}

\subsection{M\"{o}bius ladders and Petersen graph}

All M\"{o}bius ladders $M_{2m}$ are toroidal.   M\"{o}bius ladder $M_{6}=K_{3,3}$, Petersen graph  and Heawood graph 
are both, toroidal and $1$-planar.

Given the M\"{o}bius ladder $M_{2m}$, call 
{\em ring-edges} $2m$ those
belonging to the $2m$-cycle $C_{1,\dots ,2m}$,  and {\em rung-edges} all other ones, i.e., $(i,i+m)$ for $i=1,\dots ,m$.

For any odd $t$ dividing $m$, denote by $C(m,t)$ the $(m+t)$-cycle of $M_{2m}$, having, up to a cyclic shift, the form 
$$1, \dots , 1+\frac{m}{t},1+\frac{m}{t}+m, \dots , 1+\frac{2m}{t}+m,  1+\frac{2m}{t}+2m,\dots , 
1+\frac{3m}{t}+2m, \dots ,$$
i.e., $t$ consecutive sequences of $\frac{2m}{t}-1$ ring-edges, followed by a rung-edge. Such $C(m,1)$ exists for any $m\ge 3$;
 for $t> 1$, their existence requires divisibility of $m$ by $t$. Clearly, the number of $(m+t)$-cycles $C(m,t)$ is
 $\frac{2m}{t}$. 
 
 \begin{conjecture} If $G=M_{2m}$ ($m\ge 4$), then among  facets of $\CUTP(G)$ there are:
 
 two orbits of $4m$ and $2m$ edge facets (from all $2m$ ring- and $m$ rung-edges),
 
 orbit of $2^3c_4=8m$ $4$-cycles facets (from all $m$ $4$-cycles),
 
 orbit of $2^m2m$ $(m+1)$-cycle facets (from all $2m$ $(m+1)$-cycles $C(m,1)$),  
 
 for any odd divisor $t>1$ of $m$, orbit of $2^{m+t}\frac{m}{t}$ $(m+t)$-cycle facets (from all $(m+t)$-cycles $C(m,t)$).

   \end{conjecture}  
 There are no other orbits for $m=3,4$ and for $m=3$ first two orbits unite into one of $18$ edge facets, while all other orbits
 unite into  one of $2^3c_4=72$ $4$-cycle facets.  $\CUTP(M_{10})$ has only one more orbit: the orbit of $2^{10}$  facets of incidence $15$ (i.e., simplicial facets), defined by a cyclic shift of
 $$\sum_{i=1}^{10}\frac{1}{2}(3-(-1)^i)x_{i,i+1}+\sum_{i=0}^mx_{i,i+m}-2(x_{5,10}+2x_{1,2}+x_{3,8}).$$
$\CUTP(M_{12})$ also has only one more orbit: $2^{12}6$ similar facets of incidence $20$.

\bigskip

Petersen graph has three circuit double covers: 
by  six $5$-gons (actually, zigzags), by  five cycles of lengths 
$9,6,5,5,5$ and  by $5$ cycles of lengths $8,6,6,5,5$. It 
can be embedded in projective plane,  in torus 
and in Klein bottle with corresponding sets of six, five and five faces.

Petersen graph have only $5-,6-,8-$ and $9$-cycles; it has $c_5=12$ and $c_6=10$. Heawood graph, i.e., {\em $(3,6)$-cage}, have the girth $6$ and $c_6=28$, $c_8=|E|=21$. 

\begin{proposition}
$\CUTP(Petersen\, graph)$  has $3,614$ facets in $4$ orbits:
\begin{enumerate}
 
\item orbit of $2|E|=30$ edge facets, 
\item orbit of   $2^4 c_5 =192$ $5$-cycle facets, 
\item orbit of $2^5c_6 =320$  $5$-cycle facets,
\item orbit 
  of $2^{10}3$ simplexes, represented by 
  $$(C_{12345}-2x_{15}) -(C_{1'4'2'5'3'}-x_{1'4'}-x_{2'5'})  + 2\sum_{1\le i\le 5}x_{ii'},$$
   where Petersen graph is  seen as $C_{12345}+C_{1'4'2'5'3'}+ \sum_{1\le i\le 5}x_{ii'}$.   
 \end{enumerate}
 \end{proposition}  
 
\begin{remark}
Three of all $9$ orbits of facets of $\CUTP(Heawood\, graph)$,   are: 
\begin{enumerate}
\item $2|E|=42$ edge facets,
\item $2^5c_6=896$ $6$-cycle facets and
\item $2^7c_8=2,688$ $8$-cycle facets.
\end{enumerate}
 
\end{remark}

\subsection{Complete-like graphs}

$K_n$ is toroidal only for  $n=5,6,7$, while it is $1$-planar only for $n=5,6$.
   Among complete multipartite graphs $G$, the planar ones are:  $K_{2,m}$; $K_{1,1,m}$; $K_{1,2,2}$; 
$K_{1,1,1,1}=K_4$ and their subgraphs. 
The  $1$-planar $G$  are, besides above: $K_6$; $K_{1,1,1,6}$; $K_{1,1,2,3}$; $K_{2,2,2,2}$; $K_{1,1,1,2,2}$ and their subgraphs (\cite{CzapPlanarity})

Given sets $A_1,\dots , A_t$ with $t\ge 2$ and $1\le |A_1|\le \dots \le |A_t|$,  let $G$  be complete multipartite graph  $K_{a_1,\dots ,a_t}$ with $a_i=|A_i|$ for $1\le i\le t$.  


All possible chordless cycles in $G$ are $c_3=\sum_{1\le i<j<k\le t}a_ia_j a_k$.
 triangles and $c_4=\sum_{1\le i \le t}{a_i \choose 2}{a_j \choose 2 }$  quadrangles. Hence, $c_3>0$ if and only if $t>2$ and
 $c_4>0$ if and only if  $(a_1,t)\neq (1,2)$.  
  So, among edge and $s$-cycle facets of  $\CUTP(G)$, only three such orbits are possible:
 $2|E|$ edge facets if $t=2$, $4c_3$ $3$-cycle facets if $\ge 3$ and  $8c_4$ $4$-cycle facets if $(a_1,t)\neq (1,2)$.  

 All  cases, when there are no other facets, i.e., when $G$ has no $K_5$-minor, are given in Table \ref{TableSomeCUTpolytopes1}; note that the facets are simplexes for  $G =K_{2,2}$ and $K_{1,1,1,1}$.
 In particular, $G=K_{m+i}-K_m, m>1,$ has no $K_5$-minor only for $i=1,2,3$. The facets of  $\CUTP(G)$ are the orbit of $2m$ edge facets for $i=1$, the  orbit of $2m$ $3$-cycle facets for $i=2$ and two orbits (of sizes $12m$ and $4$) of $3$-cycle facets for $i=3$.   
 
Some of remaining cases presented in Table \ref{TableSomeCUTpolytopes2}. 
For $G=K_{m+4}-K_m=K_{1,1,1,1,m>1}$ and  $K_{1,1,2,m>2}$, the number of orbits stays constant
for any $m$: $4$ and $7$, respectively.   


Given sequence $b_1, \dots , b_n$
 of integers, which sum to $1$, let us call 
 $$hyp(b)=\sum_{1\le i,j \le n}x_{ij}b_ib_j\le 0$$
 (when it is applicable) {\em hypermetric inequality}. Note that $hyp(1,1,-1,0,\dots , 0)$ is usual triangle inequality.
 Denote $hyp(b)$ with all non-zero $b_i$ being $b_x=b_y=1=-b_z$ 
 by $Tr(x,y;z)$ and $hyp(b)$ with  all non-zero $b_i$ being 
 $b_x=b_y=b_z=1=-b_u=-bv$ 
 by $Pent(x,y,z;u,v)$.

If $G=K_{1,1,2,m}$ with $m\ge 3$, then 
$\CUTP(G)$ has 
 $8+20m+8{m\choose 2}(16m-15)$ facets in $7$ orbits: $3$ orbits of $8,4m, 16m$ $3$-cycle facets, one orbit of $8{m\choose 2}$ 
 $4$-cycle facets and $3$ orbits of  $64{m\choose 2},64{m\choose 2}, 384{m\choose 3}$ $\{0,\pm 1 \}$-valued non-s-cycle facets,
 having $4$ values $-1$ and $11,11,12$ values of $1$. The partition is $\{1\},\{2\},\{3,4\},\{5,\dots, m+4\}$.

$\CUTP(K_{1,1,2,2})$ has $184$ facets in $4$ orbits: 
$2$ orbits of $8+8,32$ $3$-cycle facets, one orbit of $8$ 
 $4$-cycle facets and one orbits of  $2^7$  facets, represented by
\begin{equation*}
hyp(1,1,1,-1,\,-1,0)+hyp(0,0,1,1,\,0,-0,-1)\le 0.
\end{equation*}

The graph $G=K_{m+t}-K_{m}=K_{1,\dots ,1,m}$ has a $K_5$-minor only if $t\ge 4$.
If $m\ge 3$, then 
 $\CUTP(G)$ has $2$ orbits of $4m{t\choose 2}$ and  $4{t\choose 3}$ $3$-cycle facets and, for $t<4$ only, no other facets. 
 The partition is $\{1\},\dots ,\{t\},\{t+1,\dots, t+m\}$. 

If $G=K_{m+4}-K_{m}$, then $\CUTP(G)$ has $8(8m^2-3m+2)$ facets in $4$ orbits:  $2$ orbits of $24m$,
$16$ $3$-cycle facets and $2$ orbits of sizes $16m, 128{m\choose 2}$,  represented
by $hyp(1,1,-1,-1,\,1,0,\dots ,0)\le 0, \mbox{~i.e.,~} Pent(1,2,5;3,4)$ and
\begin{equation*}
hyp(1,1,-1,0,\,1,-1,0, \dots , 0) + hyp(0,0,0,-1,\,1,1,0, \dots , 0)\le 0.
\end{equation*}

If $G=K_{m+5}-K_{m}$, then among many orbits of facets of $\CUTP(G)$,
there are $2$ orbits of $40,40m$ $3$-cycle facets and  $3$ orbits of $16, 80m, 20m(m-1)$
facets, represented, respectively, by
\begin{enumerate} 
\item $hyp(1,1,1,-1,-1,\,0,\dots , 0)\le 0$,
\item $hyp(1,1,-1,-1,0,\,1,0, \dots , 0)\le 0$
\item and $hyp(1,-1,-1,0,0,\,1,1,0, \dots, 0) +hyp(0,0,0,1,0,\,1-1,0, \dots, 0)\le 0$.
\end{enumerate}

Among $12$ remaining orbits for  $K_{7}-K_{2}$, two (of sizes $2^730, 2^760$) are $\{0,\pm 1\}$-valued;
they are represented, respectively, by
\begin{enumerate} 
\item $hyp(1,1,-1,-1,1,1,-1)+ (x_{34}+x_{47}-x_{2,7}+x_{12}-x_{13})\le 0$ and  
\item $(x_{13} + x_{34} + x_{45} + x_{15}) + (x_{23}+x_{36}+  x_{67} + x_{27}) - (x_{14} + x_{47} - x_{57} + x_{25} + x_{26} + x_{16})$.
\end{enumerate}

Let $G=Pyr^2(C_m)$. Clearly, it is $K_4, K_5$ if $m=2,3$, respectively. For $m\ge 4$, it hold $A(G)=4m$ and all chordless cycles
$3m$ triangles and unique $m$-cycle. Any of $3m+1$ edges belongs to a triangle. So, among orbits of facets of $\CUTP(G)$, there are two (of size $8m$ and $4m$) orbits of $3$-cycle facets
and  orbit of $2^{m-1}$ $m$-cycle facets. All other facets for $m\le 7$ are $\{0,\pm 1\}$-valued. 

For $Pyr^2(C_{1\dots m})$ with $m=4$, unique remaining orbit consists of $2^7$ facets, represented by
$Pent(3,5,5;1,2)+Tr(1,2;4)$. Among remaining orbits for $m=5$ and $7$, there is an orbit of $2^{m+1}$ facets represented by
\begin{enumerate} 
\item $Pyr^2(C_{12345})-2((x_{45}+x_{67})+(x_{16}+x_{17}+x_{36}+x_{37}))\le 0$ and, respectively, by
\item $Pyr^2(C_{12345})-2((x_{12}+x_{19})+(x_{29}+x_{38}+x_{49}+x_{58}x_{69}+x_{78}))\le 0$.
\end{enumerate}

For $m=5$, two remaining orbits (each of size $2^65$) are represented by 
\begin{enumerate} 
\item $C_{12345}-2x_{15}+x_{67}+((x_{17}-x_{16})-(x_{37}-x_{36})+(x_{47}-x_{46}))\le 0$ and
\item $C_{12345}-2x_{15}+x_{67}+((x_{17}-x_{16})-(x_{47}-x_{46})+(x_{57}-x_{56}))\le 0$, respectively.
\end{enumerate}

For $m=6$, one of $4$ remaining orbits (of size $2^76$) is represented by 

$C_{123456}-2x_{12}+x_{78}+((x_{17}-x_{18})+(x_{57}-x_{58})-(x_{67}-x_{68}))\le 0$.

\bigskip
Note that $K_7-C_5=Pyr^2(C_5)$. Now, $G=K_7-C_{1234}=K_{\{7\},\{6\},\{5\},\{1,3\},\{2,4\}}$ has $c_3=19$; $\CUTP(G)$ has four orbits of facets: three (of sizes $48,24,4$) of
$3$-cycle facets and one orbit of size $32$, represented by $Pent(4,5,6;2,7)$. Each of  $K_5$-minors, $K_{\{2, 4, 5,6,7\}}$ and 
$K_{\{1, 3, 5,6,7\}}$ provides $16$ of above $32$ facets.

$G=K_7-C_{7}$ has $c_3=c_4=7$; $\CUTP(G)$ has three orbits of facets: one (of size $28$) of
$3$-cycle facets, one  (of size $56$) of $4$-cycle facets and one  of size $64$, represented
by $(K_7-C_{1234567})-2(x_{15}+Path_{27364})$.

\section{Quasi-metric polytopes over graphs}\label{Sec_Quasi}

We first define the inequalities satisfied by quasi-metrics on $n$-points.

\begin{definition}\label{DEFS_Ineqs_Or_Bound}
Given a fixed $n\geq 3$ we define:

(i) The oriented triangle inequality for all $1\leq i,j,k\leq n$
\begin{equation*}
d(i,j) \leq d(i,k) + d(k,j)
\end{equation*}

(ii) The non-negativity inequality for all $1\leq i,j\leq n$ is
\begin{equation*}
d(i,j) \geq 0
\end{equation*}

(iii) A bounded oriented metric is a metric satisfying for all $1\leq i,j,k\leq n$ the inequalities
\begin{equation*}
d(j,i) + d(i,k) + d(k,j) \leq 2 \mbox{~and~} d(i,j) \leq 1.
\end{equation*}
\end{definition}

Using this we can define the cone of quasimetrics $\QMET(K_n)$ (see \cite{QuasiMetric1,QuasiMetric2}
for more details) to be the cone of oriented metrics satisfying the inequalities (i), (ii) of \ref{DEFS_Ineqs_Or_Bound}.
We define the polytope $\QMETP(K_n)$ to be
the set of metrics satisfying the inequalities of \ref{DEFS_Ineqs_Or_Bound}.

Given a subset $S\subset \{1, \dots, n\}$ we define the {\em oriented switching}:
\begin{equation*}
  F_S(d)(i,j) = \left\{\begin{array}{cl}
  1 - d(j,i) & \mbox{if~} \left\vert S\cap \{i,j\}\right\vert = 1,\\
  d(i,j)     & \mbox{otherwise}.
  \end{array}\right.
\end{equation*}
The symmetric group $Sym(n)$ acts on $\QMET(K_n)$ and define a group of size $n!$. The oriented switchings
determine and $Sym(n)$ act on $\QMETP(K_n)$ and determine a group of size $2^{n-1} n!$.

The cone $\MET(K_n)$ and polytope $\METP(K_n)$ are embedded into $\QMET(K_n)$ and $\QMETP(K_n)$ but we have
another interesting subset:

\begin{definition}\label{DEFS_Weightable}
Given $n\geq 3$ and an oriented metric $d\in \QMET(K_n)$, $d$ is called {\em weightable} if it satisfies
the following equivalent definitions:
  
(i) An oriented metric is called weightable if there exist a function $w_i$ such that for all $1\leq i,j\leq n$
\begin{equation*}
d(i,j)  + w_i = d(j,i) + w_j
\end{equation*}

(ii) For all $1\leq i,j,k\leq n$ we have
\begin{equation*}
d(i,j) + d(j,k) + d(k,i)   =   d(j,i) + d(k,j) + d(i,k)
\end{equation*}
  
\end{definition}
We thus define the cone $\WQMET(K_n)$ and polytope $\WQMETP(K_n)$ to be the set of weightable quasimetrics
of the cone $\QMET(K_n)$ and polytope $\QMETP(K_n)$. Clearly, the oriented switching preserves $\WQMETP(K_n)$.

With all those definitions we can now define the corresponding objects on graphs:

\begin{definition}
  Let $G$ be an undirected graph; we define $E(G)$ the set of edges and $Dir(E(G))$ to be the set of directed edges of $G$:

  (i) We define the cones $\QMET(G)$ and $\WQMET(G)$ to be the projections of the
  cones $\QMET(K_n)$ and $\WQMET(K_n)$ on $\mathbb{R}^{Dir(E(G))}$.

  (ii) We define the polytopes $\QMETP(G)$ and $\WQMETP(G)$ to be the projections of the
  polytopes $\QMETP(K_n)$ and $\WQMETP(K_n)$ on $\mathbb{R}^{Dir(E(G))}$.
\end{definition}

We can now give a description by inequalities of $\QMET(G)$:

\begin{theorem}\label{QMET_theorem}
For a given graph $G$ the polyhedral cone $\QMET(G)$ is defined as the set of functions $\RR^{Dir(E)}$ such that

(i) For any directed edge $e=(i,j)$ of $G$ the inequality $0\leq d(i,j)$.

(ii) For any oriented cycle $e=(v_1, v_2, \dots, v_m)$ of $G$
\begin{equation}\label{Oriented_Cycle_Ineq}
d(v_1,v_m)\leq d(v_1,v_2) + d(v_2,v_3) + \dots + d(v_{m-1},v_m)
\end{equation}

The same results holds for $\WQMET(G)$ by adding the extra condition that there exist a function $w$ such
that $d(i,j) - d(j,i) = w_i - w_j$.

\end{theorem}
\proof Our proof is adapted from the proof of \cite[Theorem 27.3.3]{DL}.
It is clear that the cycle inequalities (i) and (ii) are valid for $d\in \QMET(K_n)$
and that edges of $G$ do not occur in their expression.
Therefore, the inequalities are also valid for the projection.

The proof of sufficiency is done by induction and is more complicated. Suppose that
the result is proved for $G+e$, i.e. $G$ to which an edge $e=(i,j)$ has been added.
Suppose we have an element $x$ of $\RR^{Dir(E(G))}$ satisfying all oriented cycle
inequalities.

We need to find an antecedent of $x$, i.e. a function $y\in \RR^{Dir(E(G) + e)}$.
That is we need to find $y(i,j)$ and $y(j,i)$. 

We write $P_{i,j}$ to be the set of directed paths from $i$ to $j$ in $G$.
Assume first that $P_{i,j}\not= \emptyset$. We write
\begin{equation*}
u_{i,j} = \min_{u \in P_{i,j}} x(u)
\end{equation*}
since $x$ is non-negative, we have $u_{i,j}\geq 0$.
We then write
\begin{equation*}
l_{i,j} = \max_{v \in P_{i,j}, f\in v} x(r(f)) - x(v - f)
\end{equation*}
with $r(f)$ the reversal of the directed edge $f$.
If $P_{i,j}=\emptyset$, i.e. if the edge $e$ is connecting two connected components of $G$ then we set $l_{i,j} = u_{i,j} = 0$.

We have $l_{i,j}\leq u_{i,j}$ since otherwise we could take a path $u$ realizing the minimum $u_{i,j}$, a path $v$ and directed edge $f$ realizing the maximum $l_{i,j}$ put it together and get a counterexample to the oriented cycle inequality (ii).

So, we can find a value $y_{i,j}$ such that
\begin{equation*}
l_{i,j} \leq y_{i,j} \leq u_{i,j}
\end{equation*}
and since $u_{i,j} \geq 0$ we can choose $y_{i,j}\geq 0$. The same holds for $y_{j,i}$. Therefore we found an antecedent of $x$ in $\RR^{Dir(E(G) + e)}$ and this proves the result for $\QMET(G)$ and so the stated theorem.

For $\WQMET(G)$ we have to adjust the induction construction. If $P_{i,j} = \emptyset$ then we can adjust the values of the weights $w$ such that $w_i = w_j$. This is possible since the weights are determined up to a constant term.

On the other hand if $P_{i,j}$ is not empty then the weight is already given and we should get in the end 
$y_{i,j} - y_{j,i} = w_i - w_j$. Actually this is not a problem since it can be easily be shown that
$u_{i,j} - u_{j,i} = w_i - w_j$ and $l_{i,j} - l_{j,i} = w_i - w_j$ and so the inductive construction works. \qed

Now we turn to the construction for the polytope case. 

\begin{theorem}
For a given graph $G$ the polytope $\QMETP(G)$ is defined as the set of functions $\RR^{Dir(E)}$ such that

(i) For any directed edge $e=(i,j)$ of $G$ the inequality $0\leq d(i,j)\leq 1$ holds

(ii) For any oriented cycle $C=(v_1, v_2, \dots, v_m)$ of $G$ and subset $F$ of odd size
\begin{equation}\label{Oriented_Cycle_Ineq_Poly}
\sum_{f=(v, v') \in F} d(v', v) - \sum_{f=(v, v') \in C - F} d(v, v') \leq \vert F\vert - 1
\end{equation}

The same results holds for $\WQMETP(G)$ with the extra condition that there exist a function $w$ such
that $d(i,j) - d(j,i) = w_i - w_j$.

\end{theorem}
\proof The proof follows by remarking that the inequalities (i) and (ii) are the oriented switchings of the
non-negative inequality and oriented cycle inequality \ref{Oriented_Cycle_Ineq}. Thus the proof follow from
Theorem \ref{QMET_theorem} and the same proof strategy as \cite[Theorem 27.3.3]{DL}. \qed

The oriented multicut cones defined in the introduction are very complicated. In particular the oriented
multicuts are not stable under oriented switchings.
However, we have $\WOMCUTP(K_n) = \WQMETP(K_n)$ for $n\leq 4$.
Based on that and analogy with Theorem \ref{K5wonder} a natural conjecture would be that
$\WOMCUTP(G) = \WQMETP(G)$ if $G$ has no $K_5$ minor.
But it seems that for some other graphs with no $K_5$ minor we have $\WOMCUTP(G) \not= \WQMETP(G)$.

\section{hemi-metric polytopes over simplicial complexes}\label{Sec_Simplicial}

We can also generate metrics to a measure of distance of more than $2$ objects.
Our approach differs from \cite{DezaRosenberg1,DezaRosenberg2,DezaDutourHemimetric1,DD_mining}
and has the advantage of allowing to define it on complexes.

We consider by $Set_{n,m}$ the set of subsets of $m+1$ points of $\{1, \dots, n\}$.

\begin{definition}
  Let us fix $m\geq 1$ and $n$:

  (i) A $m$-dimensional complex is formed by a subset of $Set_{n,m}$.

  (ii) A {\em closed manifold} of dimension $m$ is formed by a subset ${\mathcal S}$ of 
  $Set_{n,m}$ such that for each subset $S$ of $m$ points of $\{1, \dots, n\}$
  the number of simplices of ${\mathcal S}$ containing $S$ is even.
\end{definition}

For the case $m=1$ the closed manifold of above definition corresponds to the closed cycles.
We now proceed to defining the corresponding cycle inequalities:

\begin{definition}
Let us fix $m\geq 1$ and $n$. Given a $m$-dimensional complex $K$ on $\{1, \dots, n\}$, the hemimetric cone $\HMET(K)$ is formed by the functions $d$ on $K$ satisfying

(i) the non-negative inequalities
\begin{equation*}
d(\Delta) \geq 0
\end{equation*}
for all $\Delta\in K$.

(ii) For all closed manifolds $(\Delta_1, \dots, \Delta_r)$ formed by simplices $\Delta_i\in K$ the inequalities
\begin{equation*}
d(\Delta_i) \leq \sum_{1\leq j\leq r, i\not= j} d(\Delta_j)
\end{equation*}
for all $1\leq i\leq r$.
\end{definition}
For $m=1$ the definition corresponds to the one of $\MET(G)$.

\begin{theorem}
Let us fix $m\geq 1$ and $n$. Let us take $K$ a $m$-dimensional complex on $n$ points.
The cone $\HMET(K)$ is the projection of $\HMET(Set_{n,m})$ on the simplices included in $K$.
\end{theorem}
\proof Our proof is adapted from the proof for metric of \cite[Theorem 27.3.3]{DL}.
The inequalities for $\HMET(K)$ are clearly valid on $HMET(Set_{n,m})$ which proves one
inclusion.

We want to prove it by induction the other inclusion.
Suppose that we have a metric $d\in \HMET(K)$ and a simplex $\Delta\notin K$.
We want to find a metric $d'$ on $\HMET(K + \Delta)$. That is we need to find
a value of $d(\Delta)$ that extends the inequality.
For a subset $S\subset Set_{n,m}$ we define 
\begin{equation*}
d(S) = \sum_{\Delta'\in S} d(\Delta').
\end{equation*}
Let us consider the 
\begin{equation*}
W_{K,\Delta} = \left\{ U\subset K \mbox{~:~} U\cup \{\Delta\} \mbox{~is~a~closed~manifold}\right\}.
\end{equation*}
We now define the upper bound
\begin{equation*}
u_{K,\Delta} = \min_{U \in W_{K, \Delta}} d(U).
\end{equation*}
We have $u_{K,\Delta}\geq 0$ since $d\in \HMET(K)$ implies $d(\Delta')\geq 0$.

The lower bound is formed by
\begin{equation*}
l_{K,\Delta} = \max_{P\in W_{K,\Delta}, F\in P} d(F) - d(P - F).
\end{equation*}
Suppose that $l_{K,\Delta} > u_{K,\Delta}$. We have $u_{K,\Delta}$ realized by $U_0$ and
$l_{K,\Delta}$ is realized by $L_0$ and a face $F_0\in L_0$. 
The union $L_0\cup U_0$ is not necessarily a closed manifold since $L_0\cup U_0$ may share
simplices.
If that is so we remove them and consider instead $W_0 = L_0 \cup U_0 - L_0 \cap U_0$.

The inequality $l_{K,\Delta} > u_{K,\Delta}$ implies then
\begin{equation*}
d(F_0) > d(L_0 - F_0) + d(U_0) = d(W_0 - F_0) + 2 d(L_0\cap U_0) \geq d(W_0 - F_0)
\end{equation*}
which violates the fact that $d\in \HMET(K)$. Thus we can find a value $\alpha$ with
\begin{equation*}
l_{K,\Delta} \leq \alpha \leq u_{K,\Delta}   \mbox{~and~} \alpha\geq 0.
\end{equation*}
Thus we can find a value for $d(\Delta)$ that is compatible with an extension. \qed

The inequality set defining $\HMET(K)$ is highly redundant but is still finite so,
the cone $\HMET(K)$ is actually polyhedral.

On the other hand, using the inequalities obtained from the simplex does not work.
Consider for example the complex $Set_{6,2}$. The Octahedron has $6$ vertices and $8$
faces and is a closed manifold. Thus it determines an inequality of the form
\begin{equation*}
x_{000} \leq x_{100} + x_{010} + x_{001} + x_{110} + x_{101} + x_{011} + x_{111}
\end{equation*}
which is not implied by the inequality on the simplices. The proof can be done by linear
programming using our software {\em polyhedral} (\cite{Polyhedral}). This proves
that our construction is different from the one of \cite{DezaRosenberg1,DezaRosenberg2,DD_mining}
and it would be interesting to redo the computations of those works.

\section{Acknowledgments}

Second author gratefully acknowledges support from the Alexander von Humboldt foundation.

\providecommand{\bysame}{\leavevmode\hbox to3em{\hrulefill}\thinspace}
\providecommand{\MR}{\relax\ifhmode\unskip\space\fi MR }
\providecommand{\MRhref}[2]{%
  \href{http://www.ams.org/mathscinet-getitem?mr=#1}{#2}
}
\providecommand{\href}[2]{#2}


\begin{thebibliography}{10}

\bibitem{ArchDeacon35}
D.~Archdeacon, \emph{A {K}uratowski theorem for the projective plane}, J. Graph
  Theory \textbf{5} (1981), no.~3, 243--246, URL:
  \url{http://dx.doi.org/10.1002/jgt.3190050305}, doi:10.1002/jgt.3190050305.

\bibitem{LeggettGargIneq}
D.~Avis, P.~Hayden, and M.~M. Wilde, \emph{Leggett-{G}arg inequalities and the
  geometry of the cut polytope}, Phys. Rev. A (3) \textbf{82} (2010), no.~3,
  030102, 4, URL: \url{http://dx.doi.org/10.1103/PhysRevA.82.030102},
  doi:10.1103/PhysRevA.82.030102.

\bibitem{OrbitFacetCutPolytope6}
D.~Avis and Mutt, \emph{All the facets of the six-point {H}amming cone},
  European J. Combin. \textbf{10} (1989), no.~4, 309--312, URL:
  \url{http://dx.doi.org/10.1016/S0195-6698(89)80002-2},
  doi:10.1016/S0195-6698(89)80002-2.

\bibitem{BarahonaCutContract}
F.~Barahona, \emph{The max-cut problem on graphs not contractible to
  {$K_{s}$}}, Oper. Res. Lett. \textbf{2} (1983), no.~3, 107--111, URL:
  \url{http://dx.doi.org/10.1016/0167-6377(83)90016-0},
  doi:10.1016/0167-6377(83)90016-0.

\bibitem{BarahonaCutsMatchingsPlanar}
F.~Barahona, \emph{On cuts and matchings in planar graphs}, Math. Programming
  \textbf{60} (1993), no.~1, Ser. A, 53--68, URL:
  \url{http://dx.doi.org/10.1007/BF01580600}, doi:10.1007/BF01580600.

\bibitem{BarahonaMahjoubOnCutPolytope}
F.~Barahona and A.~R. Mahjoub, \emph{On the cut polytope}, Math. Programming
  \textbf{36} (1986), no.~2, 157--173, URL:
  \url{http://dx.doi.org/10.1007/BF02592023}, doi:10.1007/BF02592023.

\bibitem{BT_hyperconvexity}
D.~Bryant and P.~F. Tupper, \emph{Hyperconvexity and tight-span theory for
  diversities}, Adv. Math. \textbf{231} (2012), no.~6, 3172--3198, URL:
  \url{http://dx.doi.org/10.1016/j.aim.2012.08.008},
  doi:10.1016/j.aim.2012.08.008.

\bibitem{BT_hyperconvexityOTHER}
D.~Bryant and P.~F. Tupper, \emph{Diversities and the geometry of hypergraphs},
  Discrete Math. Theor. Comput. Sci. \textbf{16} (2014), no.~2, 1--20.

\bibitem{CR}
T.~Christof and G.~Reinelt, \emph{Decomposition and parallelization techniques
  for enumerating the facets of combinatorial polytopes}, Internat. J. Comput.
  Geom. Appl. \textbf{11} (2001), no.~4, 423--437, URL:
  \url{http://dx.doi.org/10.1142/S0218195901000560},
  doi:10.1142/S0218195901000560.

\bibitem{CompendiumOptiPb}
P.~Crescenzi and V.~Kann, \emph{Approximation on the web: a compendium of {NP}
  optimization problems}, Randomization and approximation techniques in
  computer science ({B}ologna, 1997), Lecture Notes in Comput. Sci., vol. 1269,
  Springer, Berlin, 1997, pp.~111--118, URL:
  \url{http://dx.doi.org/10.1007/3-540-63248-4_10},
  doi:10.1007/3-540-63248-4\_10.

\bibitem{CzapPlanarity}
J.~Czap and D.~Hud{\'a}k, \emph{1-planarity of complete multipartite graphs},
  Discrete Appl. Math. \textbf{160} (2012), no.~4-5, 505--512, URL:
  \url{http://dx.doi.org/10.1016/j.dam.2011.11.014},
  doi:10.1016/j.dam.2011.11.014.

\bibitem{BookGenFiniteMetrCuts}
E.~Deza, M.~Deza, and M.~Dutour~Sikiri\'c, \emph{Generalizations of finite
  metrics and cuts}, World Scientific Publishing Co. Pte. Ltd., Hackensack, NJ,
  2016, URL: \url{http://dx.doi.org/10.1142/9906}, doi:10.1142/9906.

\bibitem{Relatives}
M.~Deza, V.~P. Grishukhin, and M.~Laurent, \emph{The symmetries of the cut
  polytope and of some relatives}, Applied geometry and discrete mathematics,
  DIMACS Ser. Discrete Math. Theoret. Comput. Sci., vol.~4, Amer. Math. Soc.,
  Providence, RI, 1991, pp.~205--220.

\bibitem{DezaRosenberg2}
M.-M. Deza and I.~G. Rosenberg, \emph{{$n$}-semimetrics}, European J. Combin.
  \textbf{21} (2000), no.~6, 797--806, Discrete metric spaces (Marseille,
  1998), URL: \url{http://dx.doi.org/10.1006/eujc.1999.0384},
  doi:10.1006/eujc.1999.0384.

\bibitem{DezaRosenberg1}
M.-M. Deza and I.~G. Rosenberg, \emph{Small cones of {$m$}-hemimetrics},
  Discrete Math. \textbf{291} (2005), no.~1-3, 81--97, URL:
  \url{http://dx.doi.org/10.1016/j.disc.2004.04.022},
  doi:10.1016/j.disc.2004.04.022.

\bibitem{QuasiMetric2}
M.~Deza, M.~Dutour, and E.~Panteleeva, \emph{Small cones of oriented
  semi-metrics}, Forum for {I}nterdisciplinary {M}athematics {P}roceedings on
  {S}tatistics, {C}ombinatorics \& {R}elated {A}reas ({B}ombay, 2000), vol.~22,
  2002, pp.~199--225, URL:
  \url{http://dx.doi.org/10.1080/01966324.2002.10737587},
  doi:10.1080/01966324.2002.10737587.

\bibitem{DezaDutourHemimetric1}
M.~Deza, M.~Dutour, and H.~Maehara, \emph{On volume-measure as hemi-metrics},
  Ryukyu Math. J. \textbf{17} (2004), 1--9.

\bibitem{GeneralizationHypermetric}
M.~Deza and M.~Dutour~Sikiri\'c, \emph{{The hypermetric cone and polytope on
  eight vertices and some generalizations}}, J. Symb. Comp. \textbf{to~appear}
  (2017).

\bibitem{QuasiMetric1}
M.~Deza and E.~Panteleeva, \emph{Quasi-semi-metrics, oriented multi-cuts and
  related polyhedra}, European J. Combin. \textbf{21} (2000), no.~6, 777--795,
  Discrete metric spaces (Marseille, 1998), URL:
  \url{http://dx.doi.org/10.1006/eujc.1999.0383}, doi:10.1006/eujc.1999.0383.

\bibitem{CUTsmallGraphs}
M.~Deza and M.~D. Sikiri\'c, \emph{Enumeration of the facets of cut polytopes
  over some highly symmetric graphs}, Int. Trans. Oper. Res. \textbf{23}
  (2016), no.~5, 853--860, URL: \url{http://dx.doi.org/10.1111/itor.12194},
  doi:10.1111/itor.12194.

\bibitem{DD_mining}
M.-M. Deza and M.~Dutour, \emph{Cones of metrics, hemi-metrics and
  super-metrics}, Ann. Eur. Acad. Sci. \textbf{1} (2003), 141--162.

\bibitem{newBook}
M.-M. Deza, M.~Dutour~Sikiri\'c, and M.~I. Shtogrin, \emph{Geometric structure
  of chemistry-relevant graphs}, Forum for Interdisciplinary Mathematics,
  vol.~1, Springer, New Delhi, 2015, Zigzags and central circuits, URL:
  \url{http://dx.doi.org/10.1007/978-81-322-2449-5},
  doi:10.1007/978-81-322-2449-5.

\bibitem{DL}
M.~M. Deza and M.~Laurent, \emph{Geometry of cuts and metrics}, Algorithms and
  Combinatorics, vol.~15, Springer, Heidelberg, 2010, First softcover printing
  of the 1997 original [MR1460488], URL:
  \url{http://dx.doi.org/10.1007/978-3-642-04295-9},
  doi:10.1007/978-3-642-04295-9.

\bibitem{WebPageCutPolytopes}
M.~Dutour~Sikiri\'c, \emph{Cut polytopes}, URL:
  \url{http://mathieudutour.altervista.org/CutPolytopes/}.

\bibitem{Polyhedral}
M.~Dutour~Sikiri\'c, \emph{Polyhedral}, URL:
  \url{http://mathieudutour.altervista.org/Polyhedral/}.

\bibitem{GagarinMyrvoldChambers}
A.~Gagarin, W.~Myrvold, and J.~Chambers, \emph{The obstructions for toroidal
  graphs with no {$K_{3,3}$}'s}, Discrete Math. \textbf{309} (2009), no.~11,
  3625--3631, URL: \url{http://dx.doi.org/10.1016/j.disc.2007.12.075},
  doi:10.1016/j.disc.2007.12.075.

\bibitem{GloverHuneke103}
H.~H. Glover, J.~P. Huneke, and C.~S. Wang, \emph{103 graphs that are
  irreducible for the projective plane}, J. Combin. Theory Ser. B \textbf{27}
  (1979), no.~3, 332--370, URL:
  \url{http://dx.doi.org/10.1016/0095-8956(79)90022-4},
  doi:10.1016/0095-8956(79)90022-4.

\bibitem{OrbitFacetCutPolytope7}
V.~P. Grishukhin, \emph{All facets of the cut cone {${\bf C}_n$} for {$n=7$}
  are known}, European J. Combin. \textbf{11} (1990), no.~2, 115--117, URL:
  \url{http://dx.doi.org/10.1016/S0195-6698(13)80064-9},
  doi:10.1016/S0195-6698(13)80064-9.

\bibitem{Kuratowski1930}
C.~Kuratowski, \emph{Sur le probl\`eme des courbes gauches en topologie}, Fund.
  Math. \textbf{15} (1930), 271--283.

\bibitem{SeymourMatroidMulticommodity}
P.~D. Seymour, \emph{Matroids and multicommodity flows}, European J. Combin.
  \textbf{2} (1981), no.~3, 257--290, URL:
  \url{http://dx.doi.org/10.1016/S0195-6698(81)80033-9},
  doi:10.1016/S0195-6698(81)80033-9.

\bibitem{OrbitFacetCutPolytope5}
M.~E. Tylkin (=M.~Deza), \emph{On {H}amming geometry of unitary cubes}, Soviet
  Physics. Dokl. \textbf{5} (1960), 940--943.

\bibitem{Wagner1937}
K.~Wagner, \emph{\"uber eine {E}igenschaft der ebene {K}omplexe}, Math. Annal.
  \textbf{114} (1937), 570--590.

\end{thebibliography}
\end{document}